\documentclass[compmod]{rsl_2}
\gridframe{N}

\usepackage{amsmath,amssymb,latexsym,
natbib}

\usepackage{multirow}

\usepackage[latin1]{inputenc}
\usepackage{latexsym}
\usepackage{graphicx}

\usepackage{stmaryrd}

\usepackage{hyperref}

\usepackage{cleveref}
\usepackage[justification=centering]{caption}

\usepackage{amscd}
\usepackage{mathptmx}
\usepackage{mathrsfs}
\usepackage{color}
\usepackage{xspace}
\usepackage{bussproofs}

\usepackage{tikz}

\usepackage{centernot}

\usepackage{rotating}

\allowdisplaybreaks[4]

\DeclareSymbolFont{letters}{OML}{cmm}{m}{it}

\DeclareMathAlphabet{\mathcal}{OMS}{cmsy}{m}{n}

\newcommand{\Inql}{\ensuremath{\mathbf{InqL}}\xspace}

\newcommand{\PD}{\ensuremath{\mathbf{PD}}\xspace}

\newcommand{\IPL}{\ensuremath{\mathbf{IPL}}\xspace}

\newcommand{\CPL}{\ensuremath{\mathbf{CPL}}\xspace}

\newcommand{\PT}{\ensuremath{{\mathbf{PT}_0}}\xspace}

\newcommand{\leqc}{\ensuremath{\leq_{\mathsf{c}}}\xspace}

\newcommand{\dep}{\ensuremath{=\!\!}\xspace}

\newcommand{\sor}{\ensuremath{\otimes}\xspace}
\newcommand{\bor}{\ensuremath{\vee}\xspace}

\DeclareMathOperator*{\bigsor}{\bigotimes}

\newcommand{\cnt}{\ensuremath{\divideontimes}\xspace}

\newcommand{\NN}{\ensuremath{\mathbb{N}}\xspace}

\newcommand{\logicFont}[1]{\mathcal{#1}}

\newcommand{\halfliteral}[1]{\protect\ensuremath{#1}}
\newcommand{\literal}[1]{\halfliteral{#1}\xspace}

\newcommand{\set}[3][]{\literal{\left\{#2\;\middle|\;\ifthenelse{\equal{#1}{}}{\text{#3}}{\parbox{#1}{#3}}\right\}}}

\newcommand{\logic}[1]{\literal{\logicFont{#1}}}
\newcommand{\paraLogic}[2]{\ensuremath{\logic{#1}\ifthenelse{\equal{#2}{}}{}{(#2)}}\xspace}

\newcommand{\LL}{\ensuremath{\mathsf{L}}\xspace}

\volume{0}
\issue{0}

\pyear{2015}
\pmonth{}
\doinu{10.1017/S1755020300000000}
\title[Review of Symbolic Logic]
      {Uniform Definability in Propositional Dependence Logic\footnote{The research was carried out in the Graduate School in Mathematics and its Applications of the University of Helsinki, Finland. Results of this paper were included in the dissertation  \cite{Yang_dissertation} of the author. }
}

 \author[F. Yang]
        {Fan Yang\thanks{The author would like to thank Samson Abramsky, Jouko V\"{a}\"{a}n\"{a}nen and Dag Westerst\aa hl  for insightful discussions and valuable suggestions related to this paper. 
 }}
\affil{Delft University of Technology, The Netherlands\\
fan.yang.c\@ gmail.com
}

\leftrunninghead{Uniform Definability in Propositional Dependence Logic}
\rightrunninghead{F. Yang}

\begin{document}

\maketitle

\begin{abstract}
Both propositional dependence logic and inquisitive logic are expressively complete. As a consequence, every formula with intuitionistic disjunction or intuitionistic implication can be translated equivalently into a formula in the language of propositional dependence logic without these two connectives. We show that although such a (non-compositional) translation exists, neither intuitionistic disjunction nor intuitionistic implication is uniformly definable in propositional dependence logic.
\end{abstract}


\section{Introduction}
In this paper, we study the uniform definability problem of connectives in propositional dependence logic.

\emph{Dependence logic} is a  logical formalism that characterizes the notion of ``dependence'' in social and natural sciences. 
First-order dependence logic was introduced by  \cite{Van07dl} as a development of \emph{Henkin quantifier} by \cite{henkin61} and \emph{independence-friendly logic} by  \cite{HintikkaSandu1989}. Recently, propositional dependence logic was studied and axiomatized in \cite{VY_PD,SanoVirtema2014}. With a different motivation,  \cite{InquiLog} introduced and axiomatized \emph{inquisitive logic}, which can be viewed as a variant of propositional dependence logic with intuitionistic connectives. 


Dependency relations are characterized in propositional dependence logic by a new type of atomic formulas $\dep(\vec{p},q)$, called \emph{dependence atoms}. Intuitively, the atom specifies that \emph{the truth value of $q$ is determined by those of $\vec{p}$}. 
The semantics of the logic is called \emph{team semantics}, introduced by \cite{Hodges1997a,Hodges1997b} originally as a compositional semantics for independence-friendly logic. The basic idea of this new semantics is that properties of dependence cannot be manifested in \emph{single} valuations. Therefore unlike the case of classical propositional logic, formulas in  propositional dependence logic are evaluated on \emph{sets} of valuations (called \emph{teams}) instead.

Both propositional dependence logic and inquisitive logic are expressively complete with respect to downward closed team properties, as proved in \cite{VY_PD,InquiLog}. As a consequence, every instance of the intuitionistic disjunction and the intuitionistic implication can be translated into a formula in the language of propositional dependence logic (\PD) without these connectives. In this paper, we show that although such a (non-compositional) translation exists, neither of intuitionistic disjunction and intuitionistic implication is \emph{uniformly definable} in \PD.

\newpage

 This work is inspired  by \cite{Pietro_uniform}, in which the weak universal quantifier $\forall^1$ of team semantics is shown to be not uniformly definable in first-order dependence logic, even though it is definable in the logic. Similar results are also found in  \cite{ivano_msc}, where it is proved that in inquisitive logic conjunction is definable but not uniformly definable in terms of the other connectives. Another related work  is a recent result (but in a different setting) by \cite{GorankoKuusisto16} that propositional dependence and independence logic can be translated but not compositionally translated into what the authors call propositional logics of determinacy and independence, which are logics defined on the basis of Kripke semantics instead of team semantics.

This paper is organized as follows. In Section 1, we recall the basics of propositional dependence logic and its variants. In Section 2, we define the notion of uniform definability of connectives  and discuss its connection with compositional translations between logics. 
In Section 3, we study the properties of contexts for \PD, which is a crucial notion for the main argument of the paper. Section 4 presents the main results: neither intuitionistic implication nor intuitionistic disjunction is uniformly definable in \PD. 


\section{Propositional dependence logic and its variants}\label{sec:preli}

In this section, we recall the basics of propositional dependence logic and its variants. For more details of the logics, we refer the reader to \cite{VY_PD}. 



For the purpose of this paper, let us start with recalling the following definition of the syntax of a propositional logic in general.
\begin{definition}[syntax of a propositional logic]\label{synatx_logic}
The language of a propositional logic \LL is a pair $(\textsf{Atm}_\LL,\textsf{Cnt}_\LL)$, where $\textsf{Atm}_\LL$ is a set of atoms, and $\textsf{Cnt}_\LL$ is a set of connectives (each with an arity). The set $\textsf{WFF}_\LL$ of \emph{well-formed formulas} of \LL is defined inductively as follows:
\begin{itemize}
\item   $\alpha\in\textsf{WFF}_\LL$ for all $\alpha\in \textsf{Atm}_\LL$;
\item  if $\phi_1,\dots,\phi_m\in \textsf{WFF}_\LL$ and $\cnt\in \textsf{Cnt}_\LL$ is an $m$-ary connective, then $\cnt(\phi_1,\dots,\phi_m)\in \textsf{WFF}_\LL$.
\end{itemize}
\end{definition}

In this paper, we consider \emph{propositional logics of dependence}, among which there is one, known as \emph{propositional downward closed team logic}, that has the largest set of atoms and connectives. Below we define its syntax.

\begin{definition}\label{pt_syntax_df}
Fix a set ${\rm Prop}$ of propositional variables and denote its elements by $p,q,\dots$ (possibly with subscripts).  
The language of \emph{propositional downward closed team logic} (\PT) is the pair $(\textsf{Atm}_{\PT},\textsf{Cnt}_{\PT})$, where
\begin{itemize}
\item $\textsf{Atm}_\PT=\{p,\neg p,\bot\mid p\in {\rm Prop}\}\cup\{\dep(p_{1},\dots,p_{k},q)\mid p_1,\dots,p_k,q\in {\rm Prop}\}$,
\item $\textsf{Cnt}_\PT=\{\wedge,\sor,\vee,\to\}$.
\end{itemize}
Well-formed formulas of \PT are also given by the following grammar:
\[
    \phi::= \,p\mid \neg p\mid \bot\mid\dep(p_1,\dots,p_k,q)\mid
(\phi\wedge\phi)\mid
(\phi\sor\phi)\mid(\phi\vee\phi)\mid(\phi\to\phi)
\] 
\end{definition}

We call the formulas $p,\neg p,\bot$ \emph{propositional atoms}, and the formula $\dep(p_{1},\dots,p_{k},q)$ is called a \emph{dependence atom}. 
The connectives $\sor$, $\vee$ and $\to$ are called \emph{tensor} (disjunction),  \emph{intuitionistic disjunction} and \emph{intuitionistic implication}, respectively. Note that unlike in the literature of dependence logic, where negation is usually treated as a connective that applies only to atomic formulas (i.e., formulas are assumed to be in negation normal form), for reasons that will become clear in the sequel, in this paper we view $\neg p$ as an atomic formula and do not treat negation as a connective of the logic.

Fragments of \PT formed by restricting the sets $\textsf{Atm}_\PT$ and $\textsf{Cnt}_\PT$ are called \emph{propositional logics of dependence}. We now define the languages of those propositional logics of dependence that we study in this paper.

\begin{definition}\label{syntax_pd}
The language of \emph{propositional dependence logic} (\PD) is the pair $(\textsf{Atm}_{\PD},\textsf{Cnt}_{\PD})$, where 
\[\textsf{Atm}_\PD=\textsf{Atm}_\PT\text{ and }\textsf{Cnt}_\PD=\{\wedge,\sor\}.\]

The language of \emph{inquisitive logic} (\Inql) is the pair $(\textsf{Atm}_{\Inql},\textsf{Cnt}_{\Inql})$, where 
\[\textsf{Atm}_\Inql=\{p,\bot\mid p\in {\rm Prop}\}\text{ and }\textsf{Cnt}_\Inql=\{\wedge,\vee,\to\}.\]

\end{definition}

For the semantics, propositional logics of dependence adopt \emph{team semantics}. A \emph{team} $X$ is a set of valuations, i.e., a set of functions $v:{\rm Prop}\to \{0,1\}$.  
\begin{definition}\label{TS_PD}
We  define inductively the notion of a \PT-formula $\phi$ being \emph{true} on a team $X$, denoted $X\models\phi$, as follows:
\begin{itemize}
\item $X\models p$ iff
for all $v\in X$, $v(p)=1$
  \item $X\models\neg p$  iff
for all $v\in X$, $v(p)=0$
  \item $X\models\bot$ iff $X=\emptyset$
  \item $X\models\,\dep(p_1,\dots,p_k,q)$ iff for all $v,v'\in X$ 
  \[\big[\,v(p_1)=v'(p_1),\,\dots, \,v(p_{k})=v'(p_{k})\,\big]~\Longrightarrow ~v(q)=v'(q)\]
  \item $X\models\phi\wedge\psi$ iff $X\models\phi$ and
  $X\models\psi$
  \item $X\models\phi\sor\psi$ iff there exist teams $Y,Z\subseteq X$ with $X=Y\cup Z$ such that 
  \(Y\models\phi\text{ and }Z\models\psi\)
  \item $X\models \phi\bor\psi$ iff $X\models \phi$ or $X\models\psi$
  \item $X\models \phi\to\psi$ iff for any team $Y\subseteq X$, 
  $Y\models \phi$ implies $Y\models\psi$.
\end{itemize}

\end{definition}

We write $\phi(p_1,\dots,p_n)$ to mean that the propositional variables occurring in $\phi$ are among $p_1,\dots,p_n$. If $X$ is a team and $N$ is a set of propositional variables, then we write $X\upharpoonright N=\{v\upharpoonright N\mid v\in X\}$ and call $X\upharpoonright N$ a \emph{team on $N$}.

Basic properties of \PT are listed in the theorem below; see \cite{VY_PD} for the proof.

\begin{thm}
Let $\phi(p_1,\dots,p_n)$ be a \PT-formula, and $X$ and $Y$ two teams. 
\begin{description}
\item[(Locality)] If $X\upharpoonright\{p_1,\dots,p_n\}=Y \upharpoonright\{p_1,\dots,p_n\}$, then
\(X\models\phi\iff Y\models\phi.\)
\item[(Downward Closure Property)] If $X\models\phi$ and $Y\subseteq X$, then $Y\models\phi$.
\item[(Empty Team Property)] $\emptyset\models\phi$.
\item[(Disjunction Property)] If $\models\phi\vee\psi$, then $\models\phi$ or $\models\psi$.
\end{description}
\end{thm}



For each formula $\phi(p_1,\dots,p_1)$, we write $\llbracket \phi\rrbracket_N$ for the set of all teams on $N=\{p_1,\dots,p_n\}$ that satisfies $\phi$, i.e.,
\begin{equation}\label{truth_set_N}
\llbracket \phi\rrbracket_N:=\{X\subseteq 2^N\mid X\models\phi\}.
\end{equation}
Write $\nabla_N$ for the family of all non-empty downward closed collections of teams on $N$, i.e., 
\begin{equation}\label{nabla_N}
\nabla_N=\{\mathcal{K}\subseteq 2^{\mathbf{2^N}}\mid \emptyset\in\mathcal{K},\text{ and }{X}\in \mathcal{K}\text{ and }{Y}\subseteq {X}\text{ imply }{Y}\in\mathcal{K}\}.
\end{equation}
We call a propositional logic \LL of dependence  \emph{expressively complete} with respect to downward closed team properties if for every $N=\{p_1,\dots,p_n\}$,
\[\nabla_N=\{\llbracket \phi\rrbracket_N:~\phi(p_1,\dots,p_1)\text{ is an \LL-formula}\}.\]

\begin{thm}\label{max_logic}
\PT, \PD and \Inql are expressively complete with respect to downward closed team properties.
\end{thm}
\begin{proof}
See \cite{VY_PD} and \cite{InquiLog} for the proof. We only mention here that in \cite{VY_PD} the proof of the expressive completeness of \PD makes heavy use of a formula $\Theta^\star_X$ defined for every team $X$ on $N$ and having the property that  for any team $Y$ on $N$,
\begin{equation}\label{taneli_form}
Y\models\Theta^\star_X \iff X\nsubseteq Y.
\end{equation}
If $N=\{p_1,\dots,p_n\}$ and $|X|=m+1$, then the formula $\Theta^\star_X$ is defined as
\[\Theta^\star_{X}:=\bigsor_{i=1}^m(\dep(p_1)\wedge\dots\wedge\dep(p_n))\sor\bigsor_{v\in 2^N\setminus X}(p_{1}^{v(p_1)}\wedge\dots\wedge p_{n}^{v(p_{n})}),\]
where  $\bigsor\emptyset:=\bot$. We will make use of this formula in the main argument of this paper.
\end{proof}

Sound and complete deduction systems for \PD and \Inql are defined in \cite{VY_PD,SanoVirtema2014} and \cite{InquiLog}. These systems do not admit Uniform Substitution Rule and the logics  \PD and \Inql are \emph{not} closed under uniform substitution. For instance, the \Inql-formula $\neg\neg p \to p$ is true on all teams, whereas its substitution instance $\neg\neg (p\vee\neg p)\to (p\vee\neg p)$ is not; the \PD-formula $p\sor p$ implies $p$, whereas the substitution instance $\dep(p)\sor\dep(p)$ does not imply $\dep(p)$. We will see in the sequel that the closure under uniform substitution is related to the uniform definability problem that we study in this paper.



\section{Uniform definability and compositional translations}

In this section, we define the notion of uniform definability of connectives, and discuss its connection with compositional translations between logics.

Let us start by re-examining the syntax and semantics of propositional logics in general. We defined in Definition \ref{synatx_logic} the syntax of a propositional logic \LL in general as a pair $(\textsf{Atm}_\LL,\textsf{Cnt}_\LL)$, and we also defined the syntax of propositional logics of dependence in this fashion (Definitions \ref{pt_syntax_df} and \ref{syntax_pd}). Recall also that the set of atoms of \CPL or \IPL consists of  the set {\rm Prop} (of all propositional variables) and the constant $\bot$; the set $\textsf{Cnt}_\CPL$ of connectives of \CPL contains classical negation and all the other classical connectives, and the set $\textsf{Cnt}_\IPL=\{\wedge,\vee,\to\}$ (recall:  $\neg\phi:=\phi\to\bot$). 
We now give a general definition also for the semantics of a propositional logic.

\begin{definition}[semantics of a propositional logic]
To a propositional logic \LL, we assign a class (or a set) $\nabla^\LL$ (or $\nabla$ for short) as its \emph{semantics space}. Every atom $\alpha\in \textsf{Atm}_\LL$ is associated with a set $\llbracket \alpha\rrbracket\in \nabla$, and every $m$-ary connective $\cnt\in \textsf{Cnt}_\LL$ is associated with an $m$-ary interpretation function $\pmb{\cnt}:\nabla^m\to \nabla$. The interpretation of \LL-formulas is a function $\llbracket\cdot\rrbracket^\LL:\mathsf{WWF}_\LL\to \nabla$ such that
\begin{itemize}
\item $\llbracket\alpha\rrbracket^\LL=\llbracket\alpha\rrbracket$ for every $\alpha\in \textsf{Atm}_\LL$,
\item $\llbracket\cnt(\phi_1,\dots,\phi_m)\rrbracket^\LL=\pmb{\cnt}(\llbracket\phi_1\rrbracket^\LL,\dots,\llbracket\phi_m\rrbracket^\LL)$.
\end{itemize}
\end{definition}

For a propositional logic of dependence \LL, such as \PT, \PD and \Inql, the set $\llbracket\phi\rrbracket^\LL$ consists of all of the teams that satisfies $\phi$, namely $\llbracket \phi\rrbracket^\LL:=\{X\subseteq 2^{\rm Prop}: X\models \phi\}$, and
\[
\llbracket \phi\rrbracket^\LL:=\{\mathcal{K}\subseteq 2^{\mathbf{2^{\rm Prop}}}\mid \emptyset\in\mathcal{K},\text{ and }{X}\in \mathcal{K}\text{ and }{Y}\subseteq {X}\text{ imply }{Y}\in\mathcal{K}\}.\]
Note that $\llbracket \cdot\rrbracket_N$ and $\nabla_N$ defined in equations (\ref{truth_set_N}) and (\ref{nabla_N}) in the previous section can be viewed as a restricted version of $\llbracket\cdot\rrbracket^\LL$ and $\nabla^\LL$ here in this setting.

The interpretation $\llbracket\phi\rrbracket^\CPL$ of a \CPL-formula $\phi$ is the set of all valuations that makes $\phi$ true, namely 
\(\llbracket \phi\rrbracket^\CPL:=\{v:{\rm Prop}\to 2\mid v(\phi)=1\}.\)
For an \IPL-formula $\phi$, $\llbracket\phi\rrbracket^\IPL$ is the class of all point-Kripke models that satisfies $\phi$, namely
\[\begin{split}
\llbracket \phi\rrbracket^\IPL:=\{(\mathfrak{M},w)\mid ~\mathfrak{M}\text{ is an intuitionistic Kripke model }&\text{with a node }w
 \text{ and }\mathfrak{M},w\models \phi\}.
 \end{split}\]

 A propositional logic $\LL_1=(\textsf{Atm}_{\LL_1},\textsf{Cnt}_{\LL_1})$ is said to be a \emph{sublogic} or \emph{fragment} of $\LL_2=(\textsf{Atm}_{\LL_2},\textsf{Cnt}_{\LL_2})$, written $\LL_1\subseteq \LL_2$, if  
$\textsf{Atm}_{\LL_1}\subseteq\textsf{Atm}_{\LL_2}$, $\textsf{Cnt}_{\LL_1}\subseteq \textsf{Cnt}_{\LL_2}$
and the well-formed formulas of both logics have the same interpretations in both logics (i.e., $\llbracket\phi\rrbracket^{\LL_1}=\llbracket\phi\rrbracket^{\LL_2}$ for all $\phi\in \textsf{WFF}_{\LL_1}\cap \textsf{WFF}_{\LL_2}$). In this case, if $\textsf{Atm}_{\LL_1}=\textsf{Atm}_{\LL_2}$ and $\textsf{Cnt}_{\LL_1}=\{\cnt_1,\dots,\cnt_k\}\subset \textsf{Cnt}_{\LL_1}$,
then we also write $[\cnt_1,\dots,\cnt_k]_{\LL_2}$ for $\LL_1$.

\begin{definition}\label{def_translation}
For any \LL-formulas $\phi$ and $\psi$, we write $\phi\models\psi$ if $\llbracket\phi\rrbracket^\LL\subseteq\llbracket\psi\rrbracket^\LL$. Write $\phi\equiv_\LL\psi$ (or simply $\phi\equiv\psi$) if both $\phi\models\psi$ and $\psi\models\phi$ hold. 

Let $\LL_1,\LL_2\subseteq \LL$.  The logic $\LL_1$ is said to be \emph{translatable} into $\LL_2$, in symbols $\LL_1\leq \LL_2$, if for every $\LL_1$-formula $\phi$, there exists an $\LL_2$-formula $\psi$ such that $\phi\equiv_\LL\psi$. If $\LL_1\leq \LL_2$ and $\LL_2\leq \LL_1$, then we say that $\LL_1$ and $\LL_2$ \emph{have the same expressive power}, written $\LL_1\equiv\LL_2$.
\end{definition}


Clearly, for any \PT-formulas $\phi$ and $\psi$, $\phi\equiv\psi$ iff  $X\models\phi\iff X\models\psi$ holds for all teams $X$.  An immediate consequence of Theorem \ref{max_logic} is that $\PD\equiv\Inql$, namely, \Inql and \PD are inter-translatable.

It follows from Definition \ref{def_translation} that if $\LL_1\leq\LL_2\leq \LL$, then every ($m$-ary) connective \cnt of $\LL_1$ is \emph{definable} in $\LL_2$, in the sense that for every $\LL_2$-formulas $\theta_1,\dots,\theta_m$, there exists an $\LL_2$-formula $\phi$ such that $\cnt(\theta_1,\dots,\theta_m)\equiv_\LL \phi$. We are, in this paper, however, more interested in a strengthened notion of definability of a connective, namely the \emph{uniform definability} of a connective. Closely related to this notion is a strengthened notion of translation between logics known in the literature, namely the  \emph{compositional translation} between logics (see, e.g., \cite{Rosetta94,Janssen98_comp_tran,PetersandWesterstaahl2006}).



To define \emph{uniform definability} and \emph{compositional translation} formally, let us first define the notion of \emph{context} for a logic  \LL. This definition is  inspired by that of the same notion in the first-order setting given by  \cite{Pietro_uniform}. This notion is also very similar to the notion of ``frame" by \cite{Hodges2012}; see also \cite{Hodges2016} for a comparison.

\begin{definition}[context] \label{context_def}
A \emph{context} for a propositional logic \LL is an \LL-formula with distinguished atoms $ r_i$ ($i\in \NN$). We write $\phi[r_1,\dots,r_m]$ to mean that the distinguished atoms occurring in the context $\phi$ are among $ r_1,\dots, r_m$.
For any \LL-formulas $\theta_1,\dots,\theta_m$, we write $\phi[\theta_1,\dots,\theta_m]$ for the formula $\phi(\theta_1/ r_1,\dots,\theta_m/ r_m)$. 
\end{definition}

\begin{definition}[Uniform definability of connectives]\label{uniform_defiability_def}\index{uniform definability} 
Let $\LL_1,\LL_2\subseteq \LL$. An $m$-ary connective \cnt of $\LL_1$ is said to be \emph{uniformly definable} in $\LL_2$ if there exists a context $\phi[ r_1,\dots,  r_m]$ for $\LL_2$ such that for all $\LL_2$-formulas $\theta_1,\dots,\theta_m$,
\[\phi[\theta_1,\dots,\theta_m]\equiv_\LL\cnt(\theta_1,\dots,\theta_m).\]
In this case, we say that the context  $\phi[ r_1,\dots,  r_m]$ \emph{uniformly defines} $\cnt$.
\end{definition}


The distinguished atoms $r_i$  in a context should be understood as ``place holders'' or ``holes'', which mark the places that are to be substituted uniformly by concrete instances of formulas. For the propositional logics \CPL, \IPL, \PT, \PD,  or \Inql, 
a context is a formula built from the distinguished propositional variables $r_i$ ($i\in \NN$) and other atoms using the connectives of the logic. 
For example, the formula
\[\phi[ r_1, r_2]=(\neg p_1\sor  r_1)\wedge(\dep(p_2,p_3)\sor( r_1\wedge r_2))\]
is a context for \PD. 
The formula
\(\psi[ r_1, r_2]=\neg(\neg  r_1\vee\neg r_2)\)
is a context for \CPL that uniformly defines the classical conjunction, since for any formulas $\theta_1$ and $\theta_2$ in the language of the conjunction-free fragment $[\neg,\vee]_\CPL$ of \CPL, 
\(\psi[\theta_1,\theta_2]=\neg(\neg \theta_1\vee\neg \theta_2)\equiv\theta_1\wedge\theta_2.\)

Next, we define the notion of \emph{compositional translation} in the literature using our terminology. 

\begin{definition}\label{leqc_def}\footnote{The author would like to thank Dag Westerst\aa hl for suggesting this definition. See also Section 12.2.2 in \cite{PetersandWesterstaahl2006} and \cite{Hodges2016} for similar definitions.}
Let $\LL_1,\LL_2\subseteq \LL$.  A mapping $\tau:\textsf{WFF}_{\LL_1}\to\textsf{WFF}_{\LL_2}$ is called a \emph{compositional translation} between $\LL_1$ and $\LL_2$, if the following conditions hold:
\begin{description}
\item[(i)] $\alpha\equiv_\LL\tau(\alpha)$ holds for all $\alpha\in \textsf{Atm}_{\LL_1}$;
\item[(ii)] and  for each $m$-ary connective \cnt of $\LL_1$, there is a context $\phi_\cnt[r_1,\dots,r_m]$ for $\LL_2$ which uniformly defines \cnt and
\[\tau(\cnt(\theta_1,\dots,\theta_m))=\phi_\cnt[\tau(\theta_1),\dots,\tau(\theta_m)]\]
holds for any $\LL_1$-formulas $\theta_1,\dots,\theta_m$.
\end{description}

The logic $\LL_1$ is said to be \emph{compositionally translatable} into $\LL_2$, in symbols $\LL_1\leqc\LL_2$, if there is a compositional translation $\tau$ between $\LL_1$ and $\LL_2$.
\end{definition}

The above definition implies that if $\LL_1\leqc\LL_2$, then every connective of $\LL_1$ is uniformly definable in $\LL_2$. In other words, the uniform definability of every connective of $\LL_1$ in $\LL_2$ is a necessary condition for the existence of a compositional translation from $\LL_1$ into $\LL_2$.

\begin{lem}\label{leqc2leq}
Let $\LL_1,\LL_2\subseteq \LL$. Then $\LL_1\leqc\LL_2\Longrightarrow\LL_1\leq\LL_2$. 
\end{lem}
\begin{proof}
Assume $\LL_1\leqc\LL_2$ with $\tau$ a compositional translation. It suffices to show that for each $\LL_1$-formula $\psi$,  $\psi\equiv_\LL\tau(\psi)$. 

We proceed by induction on $\psi$. If $\psi\in \textsf{Atm}_{\LL_1}$, then the required equation follows from condition (i) of the compositional translation. If $\psi=\cnt(\theta_1,\dots,\theta_m)$, where $\theta_1,\dots,\theta_m\in \textsf{WFF}_{\LL_1}$ and the context $\phi_\cnt[r_1,\dots,r_m]$ uniformly defines \cnt, then
\begin{align*}
\llbracket\cnt(\theta_1,\dots,\theta_m)\rrbracket^\LL&=\pmb{\cnt}(\llbracket\theta_1\rrbracket^\LL,\dots,\llbracket\theta_m\rrbracket^\LL)\\
&=\pmb{\cnt}(\llbracket\tau(\theta_1)\rrbracket^\LL,\dots,\llbracket\tau(\theta_m)\rrbracket^\LL)\quad(\text{by the induction hypothesis})\\
&=\llbracket\cnt(\tau(\theta_1),\dots,\tau(\theta_m))\rrbracket^\LL\\
&=\llbracket\phi_\cnt[\tau(\theta_1),\dots,\tau(\theta_m)]\rrbracket^\LL\quad\text{(since $\phi_\cnt$ uniformly defines \cnt)}\\
&=\llbracket\tau(\cnt(\theta_1,\dots,\theta_m))\rrbracket^\LL\quad(\text{since $\tau$ is a compositional translation})\\
\end{align*}
\end{proof}

However, the converse direction of Lemma \ref{leqc2leq}, i.e.,
\[``\LL_1\leq\LL_2\Longrightarrow \LL_1\leqc\LL_2 \text{''},\tag{$\ast$}\] 
 is not true in general. The next theorem by \cite{ivano_msc} is an example of the failure of ($\ast$) in propositional logics of dependence.
\begin{thm}\label{pid_non_uniform_def}
$\Inql\leq[\bot,\vee,\to]_{\Inql}$, but $\Inql\not\leqc[\bot,\vee,\to]_{\Inql}$. In particular, conjunction $\wedge$ is definable but not uniformly definable in $[\bot,\vee,\to]_{\Inql}$.  
\end{thm}
\begin{proof}
Follows from Propositions 2.5.2 and 3.5.5 in \cite{ivano_msc}.
\end{proof}

The main result of this paper is a proof that neither of intuitionistic implication $\to$ and intuitionistic disjunction $\vee$ is uniformly definable in \PD. This will then imply that $\Inql\not\leqc\PD$, even though $\Inql\leq \PD$, providing another example of the failure of ($\ast$).  

Nevertheless,  ($\ast$) does hold for most familiar logics that admit uniform substitution and have the (defined) connective $\leftrightarrow$ in the language, e.g., \CPL and \IPL. In fact, for \CPL and \IPL, the notion of a connective being definable and its being  uniformly definable coincide. A proof of this fact goes as follows: Say, for example, \cnt is a binary connective and $r_1\cnt r_2$ is equivalent to a formula $\phi(r_1,r_2,\vec{p})$ in the language of a \cnt-free sublogic $\LL_0$ of $\LL\in\{\CPL,\IPL\}$, where $\vec{p}=\langle p_1,\dots,p_n\rangle$ lists the other propositional variables involved. Then $\vdash_{\LL} (r_1\cnt r_2)\leftrightarrow \phi(r_1,r_2,\vec{p})$, which implies that $\vdash_{\LL} (\theta_1\cnt \theta_2)\leftrightarrow \phi(\theta_1,\theta_2,\vec{p})$ or $\theta_1\cnt\theta_2\equiv\phi(\theta_1,\theta_2,\vec{p})$ for any $\LL_0$-formulas $\theta_1,\theta_2$, as \LL is closed under uniform substation. From this we conclude that the context $\phi[r_1,r_2]$ for $\LL_0$ uniformly defines \cnt. It is possible to extract from the foregoing argument certain general condition under which ($\ast$) will hold. However, a propositional logic in general may have some unexpected properties that are very different from those of the familiar logics. For this reason, we leave this issue for future research and do not make any claim concerning this in this paper.

We end this section by remarking that the definitions of the familiar notion of \emph{functional completeness} of a set of connectives and \emph{independence of connectives} can be rephrased using the notion of uniform definability. A set $\{\cnt_1,\dots,\cnt_n\}$ of connectives of a propositional logic \LL is said to be  \emph{functionally complete} if and only if every connective $\cnt\in \textsf{Cnt}_\LL$ of \LL is uniformly definable in the fragment $[\cnt_1,\dots,\cnt_n]_\LL$. For example, well-known functionally complete sets of connectives of \CPL are $\{\neg,\vee\}$, $\{\neg,\wedge\}$, $\{\neg,\to\}$, $\{\mid$ (Sheffer stroke)$\}$. A connective \cnt of \LL is said to be \emph{independent} of a set $\{\cnt_1,\dots,\cnt_n\}$ of connectives of \LL, if \cnt can not be uniformly defined in the logic $[\cnt_1,\dots,\cnt_n]_\LL$. For example, all of the intuitionistic connectives $\wedge,\vee,\to$ are known to be independent of the others in \IPL.


\section{Contexts for \PD}

In this section, we investigate the properties of contexts for propositional dependence logic, which will play a crucial role in the proof of the main results of this paper.


We defined in Definition \ref{context_def} a context for \PD as a \PD-formula with distinguished propositional variables $r_i$ ($i\in \NN$) that are to be substitute uniformly by concrete instances of formulas. A subtle point that needs to be addressed here is that not only is \PD not closed under uniform substitution (as commented at the end of Section \ref{sec:preli}), but also substitution 
is not even a well-defined notion in \PD if the usual syntax is applied, since, e.g., the strings $\dep(\dep(p),q)$ and $\neg\dep(p)$ are not  well-formed formulas of \PD.  We have resolved this problem by defining a slightly different syntax for \PD (Definition \ref{syntax_pd}) than that in the literature. In particular, we do not  view negation as a connective, and dependence atoms cannot be decomposed. With this syntax, the set $\mathsf{Sub}(\phi)$ of subformulas of a context $\phi$ for \PD is defined inductively as:
\begin{itemize}
\item $\mathsf{Sub}( r_i)=\{ r_i\}$
\item $\mathsf{Sub}(p)=\{p\}$
\item $\mathsf{Sub}(\neg p)=\{\neg p\}$
\item $\mathsf{Sub}(\bot)=\{\bot\}$
\item $\mathsf{Sub}(\dep(p_1,\dots,p_k))=\{\dep(p_1,\dots,p_k)\}$
\item $\mathsf{Sub}(\psi\wedge\chi)=\mathsf{Sub}(\psi)\cup \mathsf{Sub}(\chi)\cup\{\psi\wedge\chi\}$
\item $\mathsf{Sub}(\psi\sor\chi)=\mathsf{Sub}(\psi)\cup \mathsf{Sub}(\chi)\cup\{\psi\sor\chi\}$
\end{itemize}
In this setting, a context for \PD cannot have any subformula of the form 
\[\neg r_i\text{ or }\dep(p_1,\dots,p_{m-1},r_i,p_{m+1}\dots,p_k),\] 
and thus substitution instances of a context will always be well-formed formulas of \PD.  


Two contexts $\phi[ r_1,\dots, r_{m}]$ and $\psi[ r_1',\dots, r_{m}']$ for \PD are said to be \emph{equivalent}, in symbols $\phi[ r_1,\dots, r_{m}]\approx\psi[ r_1',\dots, r_{m}']$ or simply $\phi\approx \psi$, if $ \phi[\theta_1,\dots,\theta_m]\equiv\psi[\theta_1,\dots,\theta_m]$ holds for any \PD-formulas $\theta_1,\dots,\theta_m$. A context $\phi$ is said to be \emph{inconsistent} if $\phi\approx\bot$; otherwise it is said to be \emph{consistent}.
An inconsistent context $\phi[ r_1,\dots, r_{m}]$ defines uniformly an $m$-ary connective that we shall call the \emph{contradictory connective}. The following lemma shows that we may assume that a context is either inconsistent or it does not contain a single inconsistent subformula.

\begin{lem}\label{non_bot_context_bot_free}
If $\phi[ r_1,\dots, r_{m}]$ is a consistent context for \PD, then there exists an equivalent context $\phi'[ r_1,\dots, r_{m}]$ for \PD with no single inconsistent subfromula (i.e. there is no $\psi[ r_1,\dots, r_{m}]\in \mathsf{Sub}(\phi')$ such that $\psi\approx\bot$).
\end{lem}
\begin{proof}
Assuming that $\phi[ r_1,\dots, r_{m}]$ is consistent, we find the required formula $\phi'$ by induction on $\phi$.


If $\phi[ r_1,\dots, r_{m}]$ is an atom, clearly $\phi\neq\bot$ and thus we let $\phi'=\phi$.


If $\phi[ r_1,\dots, r_{m}]=(\psi\wedge\chi)[ r_1,\dots, r_{m}]$, which is consistent, then none of $\psi$ and $\chi$ is inconsistent. By the induction hypothesis, there are contexts $\psi'[ r_1,\dots, r_{m}]$ and $\chi'[ r_1,\dots, r_{m}]$ such that $\psi'\approx\psi$, $\chi'\approx\chi$
and none of $\psi'$ and $\chi'$ contains a single inconsistent formula. Let $\phi'[ r_1,\dots, r_{m}]=\psi'\wedge\chi'$. Clearly, $(\psi\wedge\chi)\approx(\psi'\wedge\chi')$ and $(\psi'\wedge\chi')\not\approx\bot$ (for $(\psi\wedge\chi)\not\approx\bot$). Thus, by the induction hypothesis, the set $\mathsf{Sub}(\psi'\wedge\chi')=\mathsf{Sub}(\psi')\cup \mathsf{Sub}(\chi')\cup\{\psi'\wedge\chi'\}$
does not contain a single inconsistent element.


If $\phi[ r_1,\dots, r_{m}]=(\psi\sor\chi)[ r_1,\dots, r_{m}]$, which is consistent, then $\psi$
 and $\chi$ cannot be both inconsistent.  There are the following two cases:

Case 1: Only one of $\psi$ and $\chi$ is inconsistent. Without loss of generality, we may assume that $\psi[ r_1,\dots, r_{m}]\approx\bot$ and $\chi[ r_1,\dots, r_{m}]\approx \chi'[ r_1,\dots, r_{m}]$ for some context $\chi'$ for \PD that does not contain a single inconsistent subformula. Clearly, $(\psi\sor\chi)\approx (\bot\sor\chi')\approx\chi'$.
Thus, we may let $\phi'=\chi'$. 

Case 2: $\psi[ r_1,\dots, r_{m}]\approx \psi'[ r_1,\dots, r_{m}]$ and $\chi[ r_1,\dots, r_{m}]\approx \chi'[ r_1,\dots, r_{m}]$ for some contexts $\psi'$ and $\chi'$ for \PD that do not contain a single inconsistent subformula. Let $\phi'=\psi'\sor\chi'$. Clearly, $(\psi\sor\chi)\approx(\psi'\sor\chi')$ and $(\psi'\sor\chi')\not\approx\bot$ (for $(\psi\sor\chi)\not\approx\bot$).
Thus, by the induction hypothesis, the set $\mathsf{Sub}(\psi'\sor\chi')=\mathsf{Sub}(\psi')\cup \mathsf{Sub}(\chi')\cup\{\psi'\sor\chi'\}$
does not contain a single inconsistent element.
\end{proof}

Contexts for \PD are \emph{monotone} in the sense of the following lemma.

\begin{lem}\label{context_pd_monotone}
Let $\phi[ r_1,\dots, r_{m}]$ be a context for \PD and $\theta_1,\dots,\theta_m,\theta_1',\dots,\theta_m'$ \PD-formulas. If $\theta_i\models\theta_i'$ for all $1\leq i\leq m$, then $\phi[\theta_1,\dots,\theta_m]\models \phi[\theta_1',\dots,\theta_m']$.
\end{lem}
\begin{proof}
Suppose $\theta_i\models\theta_i'$ for all $1\leq i\leq m$. We prove that $\phi[\theta_1,\dots,\theta_m]\models \phi[\theta_1',\dots,\theta_m']$ by induction on $\phi$.

For the only interesting case $\phi[ r_1,\dots, r_{m}]= r_i$ ($1\leq i\leq m$), if $X\models r_i[\theta_1,\dots,\theta_m]$ for some team $X$, then $X\models\theta_i\models\theta_i'$. Thus, $X\models  r_i[\theta_1',\dots,\theta_m']$.
\end{proof}

\begin{cor}\label{context_top_imply}
For any consistent context $\phi[ r_1,\dots, r_{m}]$ for \PD, there exists a non-empty team $X$ such that $X\models\phi[\top,\dots,\top]$.
\end{cor}
\begin{proof}
Since $\phi[ r_1,\dots, r_{m}]\not\approx\bot$, there exist formulas $\theta_1,\dots,\theta_m$ and a non-empty team $X$ such that $X\models\phi[\theta_1,\dots,\theta_m]$.
As $\theta_i\models \top$ for all $1\leq i\leq  m$, by Lemma \ref{context_pd_monotone}, we obtain that $X\models\phi[\top,\dots,\top]$.
\end{proof}

In the main proofs of this paper, we will make use of the syntax trees of contexts for \PD. We assume that the reader is familiar with the notion of a syntax tree of a formula and will therefore only recall  its informal definition.

The syntax tree of a \PD-formula $\phi$ is a quadruple $\mathfrak{T}_\phi=(T,\prec ,w,\mathsf{f})$ (see Figure \ref{syntext_tree_fig} for an example) such that $(T,\prec,w)$ is a (finite) full binary tree with root $w$ (i.e., a tree in which every node has either $0$ or $2$ children) and $\mathsf{f}:T\to \mathsf{Sub}(\phi)$ is a labelling function satisfying the following conditions:
\begin{description}
\item[(i)] $\mathsf{f}(w)=\phi$
\item[(ii)] $\mathsf{f}[T]=\mathsf{Sub}(\phi)$
\item[(iii)] If $x$ is a node with $\mathsf{f}(x)=\psi\cnt\chi$ ($\cnt\in\{\wedge,\sor\}$), then $x$ has two children $y$ and $z$ with $\mathsf{f}(y)=\psi$ and $\mathsf{f}(z)=\chi$.
\end{description}
As usual, we call a node $x$ an \emph{ancestor} of a node $y$ if $x\prec y$. 
The \emph{depth} $\mathsf{d}(x)$ of a node $x$ in the tree is defined inductively as follows: $\mathsf{d}(w)=0$; if $y$ is a child of $x$, then $\mathsf{d}(y)=\mathsf{d}(x)+1$.

The value $\mathsf{f}(x)\in \mathsf{Sub}(\phi)$ of a node $x$ is also called the \emph{label} of the node $x$ or the \emph{formula associated with} $x$. Clearly, the leaf nodes (i.e., nodes with no children) are always labelled with atoms, and the labelling function $\mathsf{f}$ is in general \emph{not} one-to-one since the same subformula $\psi$ may have more than one occurrences in a formula $\phi$.

\begin{figure}[t]
\begin{center}
\begin{tikzpicture}

\node (0f)  at (0, -0.2) {$(\neg p_1\sor  r_1)\wedge\big(\!\dep(p_2,p_3)\sor( r_1\wedge r_2)\big)$};

\node (1f)  at (-3.2, -1.6) {{\footnotesize$\neg p_1\sor  r_1$}};
\node (2f)  at (-3.8, -3.7) {{\footnotesize$\neg p_1$}};
\node (3f)  at (-1.2, -3.7) {{\footnotesize$ r_1$}};
\node (4f)  at (0.7, -3.7) {{\footnotesize$\dep(p_2,p_3)$}};
\node (8f)  at (4, -1.6) {{\footnotesize$\!\dep(p_2,p_3)\sor( r_1\wedge r_2)$}};
\node (7f)  at (4.2, -3.7) {{\footnotesize$r_1\wedge r_2$}};
\node (5f)  at (2.1, -5.8) {{\footnotesize$ r_1$}};
\node (6f)  at (4.9, -5.8) {{\footnotesize$ r_2$}};

\node(0t) at (0,-.8) {$w$};
\node[draw, circle, fill, scale=0.5] (1t) at (-2.5,-2) {};
\node[draw, circle, fill, scale=0.5] (8t) at (2.5,-2) {};
\node[draw, circle, fill, scale=0.5] (2t) at (-3.5,-4) {};
\node[draw, circle, fill, scale=0.5] (3t) at (-1.5,-4) {};
\node[draw, circle, fill, scale=0.5] (4t) at (1.5,-4) {};
\node[draw, circle, fill, scale=0.5] (7t) at (3.5,-4) {};
\node[draw, circle, fill, scale=0.5] (5t) at (2.5,-6) {};
\node[draw, circle, fill, scale=0.5] (6t) at (4.5,-6) {};

\draw[line width=1pt] (0t) -- (1t);
\draw[line width=1pt] (0t) -- (8t);

\draw[line width=1pt] (1t) -- (2t);
\draw[line width=1pt] (1t) -- (3t);

\draw[line width=1pt] (8t) -- (4t);
\draw[line width=1pt] (8t) -- (7t);

\draw[line width=1pt] (7t) -- (5t);
\draw[line width=1pt] (7t) -- (6t);

\end{tikzpicture}
\caption{The syntax tree of context $(\neg p_1\sor  r_1)\wedge(\!\dep(p_2,p_3)\sor( r_1\wedge r_2))$}
\label{syntext_tree_fig}
\end{center}
\end{figure}


If $X\models \phi[\theta_1,\dots,\theta_m]$, then each occurrence of a subformula of $\phi[\theta_1,\dots,\theta_m]$ is satisfied by a subteam of $X$. This can be described explicitly by a function $\tau$, called \emph{truth function}, which maps each node in the syntax tree $\mathfrak{T}_\phi$ to a subteam of $X$ satisfying the formula associated with the node. 

\begin{definition}[Truth Function]\label{truth_function_def}
Let $\phi[ r_1,\dots, r_{m}]$ be a context for \PD with the syntax tree $\mathfrak{T}_\phi=(T,\prec ,w,\mathsf{f})$, and $\theta_1,\dots,$\allowbreak$\theta_m$ \PD-formulas. Let $N$  be the set of all propositional variables occurring in the formula $\phi[\theta_1,\dots,\theta_m]$. A function $\tau:T\to \wp(2^N)$ is called a \emph{truth function} for $\phi[\theta_1,\dots,\theta_m]$ if the following conditions hold:
\begin{description}
\item[(i)] $\tau(x)\models \mathsf{f}(x)[\theta_1,\dots,\theta_m]$ for all $x\in T$;
\item[(ii)] if $\mathsf{f}(x)=\psi\wedge\chi$ and $y, z$ are the two children of $x$, then 
\(\tau(x)=\tau(y)=\tau(z);\)
\item[(iii)] if $\mathsf{f}(x)=\psi\sor\chi$ and $y, z$ are the two children of $x$, then 
\(\tau(x)=\tau(y)\cup\tau(z).\)
\end{description}
A truth function $\tau$ such that $\tau(w)=X$ is called a \emph{truth function  over $X$}.
\end{definition}

\begin{lem}\label{truth_function_ancestor_monotone}
Let $\tau$ be a truth function for $\phi[\theta_1,\dots,\theta_m]$. If $x,y$ are two nodes in the syntax tree $\mathfrak{T}_\phi$ with $x\prec y$, then $\tau(y)\subseteq \tau(x)$. In particular, if $\tau$ is a truth function  over a team $X$, then for all nodes $x$ in $\mathfrak{T}_\phi$, $\tau(x)\subseteq X$.
\end{lem}
\begin{proof}
A routine proof by induction on $\mathsf{d}(y)-\mathsf{d}(x)$.
\end{proof}

First-order dependence logic has a game-theoretic semantics with perfect information games played with respect to teams (see Section 5.2 in \cite{Van07dl}). With obvious adaptions, one can define a game-theoretic semantics for propositional dependence logic.\footnote{In Definition 5.10 in \cite{Van07dl}, leave out game rules for quantifiers and make obvious modifications to the game rules for atoms.} A truth function defined in Definition \ref{truth_function_def} corresponds to a \emph{winning strategy} for the \textsf{Verifier} in the game. An appropriate semantic game for \PD has the property that $X\models\phi$ if and only if the \textsf{Verifier} has a winning strategy in the corresponding game. The next theorem states essentially the same property for truth functions. Cf.  Lemma 5.12, Proposition 5.11 and Theorem 5.8 in  \cite{Van07dl}.

\begin{thm}\label{PD_truth_function_context}
Let $\phi[ r_1,\dots, r_{m}]$ be a context for \PD, $\theta_1,\dots,\theta_m$ \PD-formulas and $N$ the set of all propositional variables occurring in the formula $\phi[\theta_1,\dots,\theta_m]$.
For any team ${X}$ on $N$, ${X}\models\phi[\theta_1,\dots,\theta_m]$ iff there exists a truth function $\tau$ for $\phi[\theta_1,\dots,\theta_m]$ over ${X}$.
\end{thm}
\begin{proof}
The direction ``$\Longleftarrow$'' follows readily from the definition. For the other direction ``$\Longrightarrow$'', suppose ${X}\models\phi[\theta_1,\dots,\theta_m]$. Let $\mathfrak{T}_\phi=(T,\prec ,w,\mathsf{f})$ be the syntax tree of $\phi$. We define the value of $\tau$ on each node $x$ of $\mathfrak{T}_\phi$ and verify conditions (i)-(iii) of Definition \ref{truth_function_def} by induction on the depth of the nodes.

 If $x$ the root, then define $\tau(x)={X}$. Since ${X}\models\phi[\theta_1,\dots,\theta_m]$, condition (i) is satisfied for the root $x$.

 Suppose $x$ is not a leaf node, $\tau(x)$ has been defined already and conditions (i)-(iii) are satisfied for $x$. Let $y,z$ be the two children of $x$ with $\mathsf{f}(y)=\psi$ and $\mathsf{f}(z)=\chi$ for some subformulas $\psi,\chi$ of $\phi$. We distinguish two cases.
 
Case 1: $\mathsf{f}(x)=\psi\wedge\chi$. Define 
\(\tau(y)=\tau(z)=\tau(x).\) 
Then condition (ii) for $y,z$ is satisfied. By the induction hypothesis, 
\(\tau(x)\models(\psi\wedge\chi)[\theta_1,\dots,\theta_m].\)
Thus 
 $\tau(y)\models\psi[\theta_1,\dots,\theta_m]$ and $\tau(z)\models\chi[\theta_1,\dots,\theta_m]$,
 namely condition (i) is satisfied for $y,z$.

Case 2: $\mathsf{f}(x)=\psi\sor\chi$. By the induction hypothesis, 
\(\tau(x)\models(\psi\sor\chi)[\theta_1,\dots,\theta_m].\)
Thus there exist teams $Y,Z\subseteq \tau(x)$ on $N$ such that $\tau(x)={Y}\cup {Z}$,
 ${Y}\models\psi[\theta_1,\dots,\theta_m]$ and ${Z}\models\chi[\theta_1,\dots,\theta_m]$.
Define $\tau(y)={Y}$ and $\tau(z)={Z}$. Then, conditions (i) and (ii) for $y,z$ are satisfied.
\end{proof}

The next lemma shows that a truth function is determined by its values on the leaves of the syntax tree.

\begin{lem}\label{truth_function_determine_Xi}
Let $\phi[ r_1,\dots, r_{m}]$ be a context for \PD with the syntax tree $\mathfrak{T}_\phi=\mathop{(T,\prec,w,\mathsf{f})}$, $\theta_1,\dots,\theta_m$ \PD-formulas and $N$ the set of all propositional variables occurring in the formula $\phi[\theta_1,\dots,\theta_m]$.
 If $\tau:T\to\wp(2^N)$ is a function satisfying conditions (ii) and (iii) in Definition \ref{truth_function_def} and condition (i) with respect to $\theta_1,\dots,\theta_m$ for all leaf nodes, then $\tau$ is a truth function for $\phi[\theta_1,\dots,\theta_m]$.
\end{lem}
\begin{proof}
It suffices to prove that $\tau$ satisfies condition (i) with respect to $\theta_1,\dots,\theta_m$ for all nodes $x$ of $\mathfrak{T}_\phi$. We show this by induction on the depth of $x$.

Leaf nodes satisfy condition (i) by the assumption. Now, assume that $x$ is not a leaf. Then $x$ has two children $y,z$ with $\mathsf{f}(y)=\psi$ and $\mathsf{f}(z)=\chi$ for some subformulas $\psi,\chi$ of $\phi$. Since $\mathsf{d}(y),\mathsf{d}(z)>\mathsf{d}(x)$, by the induction hypothesis, we have 
\begin{equation}\label{truth_function_determine_Xi_eq1}
\tau(y)\models \psi[\theta_1,\dots,\theta_m]\text{ and }\tau(z)\models \chi[\theta_1,\dots,\theta_m].
\end{equation}
If $\mathsf{f}(x)=\psi\wedge\chi$, then by condition (ii), $\tau(x)=\tau(y)=\tau(z)$, and $\tau(x)\models (\psi\wedge\chi)[\theta_1,\dots,\theta_m]$ follows from (\ref{truth_function_determine_Xi_eq1}). If $\mathsf{f}(x)=\psi\sor\chi$, then by condition (iii), $\tau(x)=\tau(y)\cup\tau(z)$, and $\tau(x)\models (\psi\sor\chi)[\theta_1,\dots,\theta_m]$ follows again from (\ref{truth_function_determine_Xi_eq1}).
%
%
%
\end{proof}

\section{Non-uniformly definable connectives in \PD}

In this section, we prove that neither intuitionistic implication nor intuitionistic disjunction is uniformly definable in \PD. 


Contexts for \PD are monotone (by Lemma \ref{context_pd_monotone}), thus \PD cannot define uniformly non-monotone connectives. Below we show that intuitionistic implication is not uniformly definable in \PD as it is not monotone.\footnote{The author would like to thank Samson Abramsky for pointing out this proof idea.}
\begin{thm}\label{imp_non_uniform_def}
Intuitionistic implication is not uniformly definable in \PD.
\end{thm}
\begin{proof}
Suppose that there was a context $\phi[ r_1, r_2]$ for \PD  which defines the intuitionistic implication uniformly. Then for any \PD-formulas $\psi$ and $\chi$,
\begin{equation}\label{uniform_imp_eq1}
\phi[\psi,\chi]\equiv\psi\to\chi.
\end{equation}
Clearly
\(X\models\bot\to \bot\text{ and }X\not\models \top\to\bot\)
hold for any non-empty team $X$. It follows from (\ref{uniform_imp_eq1}) that
\(X\models\phi[\bot,\bot]\text{ and }X\not\models \phi[\top,\bot].\)
But this contradicts Lemma \ref{context_pd_monotone} as $\bot\models\top$.
\end{proof}

We now proceed to give another sufficient condition for a connective being not uniformly definable in \PD, from which it will follow  that  intuitionistic disjunction is not uniformly definable in \PD. We start with a simple lemma whose proof is left to the reader.

\begin{lem}\label{sor_depth_fact1}
Let $\phi[ r_1,\dots, r_m]$ be a context for \PD and $\theta_1,\dots,\theta_m$ \PD-formulas. Let $\tau$ be a truth function for $\phi[\theta_1,\dots,\theta_m]$ over a team $X$. In the syntax tree $\mathfrak{T}_\phi$ of $\phi$, if a node $x$ has no ancestor node with a label of the form $\psi\sor\chi$, then $\tau(x)={X}$.
\end{lem}
\begin{proof}
Easy, by induction on the depth of $x$.
\end{proof}

Since, e.g., $\bot\vee\top\not\models\bot$ and $\top\vee\bot\not\models\bot$, from the above lemma it  follows that in the syntax tree of a context $\phi[r_1,r_2]$ for \PD that defines $\vee$ (if exists) every leaf node labeled with $r_1$ or $r_2$ must have an ancestor node labeled with $\sor$. Below we prove this observation  in a more general setting.

\begin{lem}\label{uniform_def_cond1}
Let \cnt be an $m$-ary connective such that for every $1\leq i\leq  m$, there are some \PD-formulas $\theta_1,\dots,\theta_m$ satisfying
\begin{equation}\label{uniform_def_cond1_eq1}
\cnt(\theta_1,\dots,\theta_m)\not\models\theta_i.
 \end{equation}
If $\phi[ r_1,\dots, r_m]$ is a context for \PD which uniformly defines \cnt, then in the syntax tree $\mathfrak{T}_\phi=(T,\prec,w,\mathsf{f})$, every leaf node labeled with $ r_i$ ($1\leq i\leq  m$) has an ancestor node with a label of the form $\psi\sor\chi$. 
\end{lem}
\begin{proof}
Suppose there exists a leaf node $x$ labeled with $ r_i$ which has no ancestor node with a label of the form $\psi\sor\chi$. By assumption, there exist \PD-formulas $\theta_1,\dots,\theta_m$ satisfying (\ref{uniform_def_cond1_eq1}) for $i$. Let $N$  be the set of all  propositional variables occurring in the formula $\phi[\theta_1,\dots,\theta_m]$. Take a team $X$ on $N$ such that
\(X\models \cnt(\theta_1,\dots,\theta_m)\text{ and }X\not\models\theta_i.\)
Since $\phi[ r_1,\dots, r_m]$ uniformly defines \cnt, we have $\cnt(\theta_1,\dots,\theta_m)\equiv\phi[\theta_1,\dots,\theta_m]$, 
implying $X\models \phi[\theta_1,\dots,\theta_m]$. By Theorem \ref{PD_truth_function_context}, there is a truth function $\tau$ for \allowbreak$\phi[\theta_1,\dots,\theta_m]$ over $X$. By the property of $x$ and Lemma \ref{sor_depth_fact1}, $\tau(x)=X$. Thus
\(X\models  r_i[\theta_1,\dots,\theta_m]\), i.e., $X\models\theta_i$; 
 a contradiction.
\end{proof}

The following elementary set-theoretic lemma will be used in the proof of Lemma \ref{PD_context_transition_reduce_lm}.

\begin{lem}\label{proper_split_lm}
Let $X,Y,Z$ be sets such that $|X|>1$, $Y,Z\neq \emptyset$ and $X=Y\cup Z$. Then there exist $Y',Z'\subsetneq X$ such that $Y'\subseteq Y$, $Z'\subseteq Z$ and $X=Y'\cup Z'$.
\end{lem}
\begin{proof}
If $Y,Z\subsetneq X$, then taking $Y'=Y$ and $Z'=Z$, the lemma holds. Now, assume that one of $Y,Z$ equals $X$. 

Case 1: $Y=Z=X$. Pick an arbitrary $a\in X$. Let $Y'=X\setminus\{a\}\subsetneq X$ and $Z'=\{a\}$. Since $|X|>1$, we have that $Z'\subsetneq X$. Clearly, $X=(X\setminus\{a\})\cup \{a\}$.

Case 2: Only one of $Y$ and $Z$ equals $X$. Without loss of generality, we assume that $Y=X$ and $Z\subsetneq X$. Let $Y'=X\setminus Z$ and $Z'=Z$. Clearly, $X=(X\setminus Z)\cup Z$ and $Y',Z'\subsetneq X$, as $\emptyset\neq Z\subsetneq X$.
\end{proof}

Next, we prove a crucial technical lemma for the main theorem (Theorem \ref{non_uniform_def}) of this section.

\begin{lem}\label{PD_context_transition_reduce_lm}
Let $\phi[ r_1,\dots, r_m]$ be a consistent context for \PD such that in the syntax tree $\mathfrak{T}_\phi=(T,\prec,w,\mathsf{f})$ of $\phi$, every leaf node labeled with $ r_i$ ($1\leq i\leq  m$) has an ancestor node labeled with a formula of the form $\psi\sor\chi$. Let $N$ be the set of all propositional variables occurring in the formula $\phi[\top,\dots,\top]$. If  $2^N\models\phi[\top,\dots, \top]$,
then there exists a truth function $\tau$ for $\phi[\top,\dots, \top]$ over $2^N$ such that $\tau(x)\subsetneq 2^N$
for all leaf nodes $x$ labeled with $ r_i$ ($1\leq i\leq  m$).
\end{lem}
\begin{proof}
By Lemma \ref{non_bot_context_bot_free}, we may assume that $\phi[ r_1,\dots, r_m]$ does not contain a single inconsistent subformula. Suppose $2^N\models\phi[\top,\dots, \top]$. The required truth function $\tau$ over $2^N$ is defined inductively on the depth of the nodes in the syntax tree $\mathfrak{T}_\phi$ in the same way as in the proof of Theorem \ref{PD_truth_function_context}, except for the following case.

For each leaf node labeled with $ r_i$, consider its ancestor node $x$ with $\mathsf{f}(x)=(\psi\sor\chi)$ of minimal depth, where $\psi,\chi\in \mathsf{Sub}(\phi)$ (the existence of such $x$ is guaranteed by the assumption). Let $y,z$ be the two children of $x$. Assuming that $\tau(x)$ has been defined already, we now define $\tau(y)$ and $\tau(z)$.

By the induction hypothesis, 
\(\tau(x)\models(\psi\sor\chi)[\top,\dots, \top].\)
The minimality of $x$ implies that $x$ has no ancestor node labeled with $\theta_0\sor\theta_1$ for some $\theta_0,\theta_1$. Thus $\tau(x)=2^N$ by Lemma \ref{sor_depth_fact1}, and there exist teams $Y_0,Z_0\subseteq \tau(x)=2^N$ such that $2^N={Y_0}\cup {Z_0}$,
 \({Y_0}\models\psi[\top,\dots, \top]\text{ and }{Z_0}\models\chi[\top,\dots, \top].\) 

\vspace{6pt}

\noindent Claim: There are non-empty teams $Y,Z$ such that  $2^N=Y\cup Z$ and
\begin{equation}\label{PD_context_transition_reduce_lm_qe1}
Y\models \psi[\top,\dots, \top]\text{ and }Z\models \chi[\top,\dots, \top].
\end{equation}

\noindent\emph{Proof of Claim.} If $Y_0,Z_0\neq\emptyset$, then taking $Y=Y_0$ and $Z=Z_0$ the claim holds. Now, suppose one of $Y_0,Z_0$ is empty. Without loss of generality, we may assume that $Y_0=\emptyset$. Then let $Z:=Z_0=2^N$. Since $\psi[ r_1,\dots, r_m]\not\approx\bot$, by Corollary \ref{context_top_imply} and the locality of \PD, there exists a non-empty team $Y\subseteq 2^N$ such that 
$Y\models \psi[\top,\dots,\top]$, as required.\hfill$\dashv$

\vspace{6pt}
Now, since $|2^N|>1$, by Lemma \ref{proper_split_lm}, there are teams $Y',Z'\subsetneq 2^N$ such that $Y'\subseteq Y$, $Z'\subseteq Z$ and $Y'\cup Z'=2^N$. Define $\tau(y)=Y'$ and $\tau(z)=Z'$.
Clearly, condition (iii) of Definition \ref{truth_function_def} for $y,z$ is satisfied. Moreover, by the downwards closure property, it follows from (\ref{PD_context_transition_reduce_lm_qe1}) that condition (i) for $y$ and $z$ is also satisfied. Hence, such defined $\tau$ is a truth function for $\phi[\top,\dots,\top]$ over $2^N$.

\vspace{6pt}
It remains to check that $\tau(x)\subsetneq 2^N$ for all leaf nodes $x$ labeled with $ r_i$ ($1\leq i\leq  m$). By the assumption, there exists an ancestor $y$ of $x$ labeled with $(\psi\sor\chi)$ of minimal depth. One of $y$'s two children, denoted by $z$, must be an ancestor of $x$ or $z=x$. Thus, by Lemma \ref{truth_function_ancestor_monotone} and the construction of $\tau$, we obtain that $\tau(x)\subseteq \tau(z)\subsetneq 2^N$.
\end{proof}

Now, we give the intended sufficient condition for a non-contradictory connective being not uniformly definable in \PD. In the proof, we will make use of the formula $\Theta^\star_X$ from the proof of Theorem \ref{max_logic} which  has the property (\ref{taneli_form}). The conditions in the statement of the next theorem are all generalized from the corresponding properties of intuitionistic disjunction, which are given in the proof of Theorem \ref{bor_non_uniform_def}. The reader is recommended to consult the proof of Theorem \ref{bor_non_uniform_def} for a better understanding of the conditions.


\begin{thm}\label{non_uniform_def}
Every non-contradictory  $m$-ary connective \cnt satisfying the following conditions  is not uniformly definable in \PD:
\begin{description}
\item[(i)] For every $1\leq i\leq  m$, there exist \PD-formulas $\theta_1,\dots,\theta_m$ such that
$\cnt(\theta_1,\dots,\theta_m)\not\models \theta_i$.
\item[(ii)] There are \PD-formulas $\delta_1,\dots,\delta_m$ such that $\models\cnt(\delta_1,\dots,\delta_m)$. 
\item[(iii)] For any finite set $N$ of propositional variables, there exist $1\leq j_1<\dots<j_k\leq m$ such that
 \begin{equation}\label{non_uniform_def_eq4}
 2^N\not\models \cnt(\alpha_1,\dots,\alpha_m),
 \end{equation}
 and for each $1\leq i\leq  m$,
\begin{equation}\label{non_uniform_def_eq3}
 \alpha_{i}=\begin{cases}
 \Theta_{2^N}^\star,&\text{ if }i=j_a~,1\leq a\leq k\\
 \top,&\text{ otherwise.}
 \end{cases}
 \end{equation}
\end{description}
\end{thm}
\begin{proof}
Suppose that $\cnt$ was uniformly definable in $\PD$ by a context $\phi[ r_1,\dots, r_m]$ for \PD such that for all \PD-formulas $\theta_1,\dots,\theta_m$, 
\begin{equation}\label{non_uniform_def_eq2}
\phi[\theta_1,\dots,\theta_m]\equiv\cnt(\theta_1,\dots,\theta_m).
\end{equation}
Since $\cnt$ satisfies condition (i), by Lemma \ref{uniform_def_cond1}, in the syntax tree $\mathfrak{T}_\phi=(T,<,r,\mathsf{f})$ of $\phi[ r_1,\dots, r_m]$, each node labeled with $ r_i$ ($1\leq i\leq  m$) has an ancestor node labeled with a formula of the form $\psi\sor\chi$.

Let $\delta_1,\dots,\delta_m$ be the \PD-formulas with $\models\cnt(\delta_1,\dots,\delta_m)$ as given by condition (ii).
By (\ref{non_uniform_def_eq2}), we have
\(\models\phi[\delta_1,\dots,\delta_m]\).
Since $\delta_i\models\top$ for all $1\leq i\leq  m$, by Lemma \ref{context_pd_monotone},
\(\models\phi[\top,\dots,\top].\)
Let $N$ be the set of all propositional variables occurring in $\phi[\top,\dots,\top]$. We have that
\(2^N\models\phi[\top,\dots,\top].\)
Since \cnt is a non-contradictory connective, $\phi[ r_1,\dots, r_m]$ is a consistent context. Then, by Lemma \ref{PD_context_transition_reduce_lm}, there exists a truth function $\tau$ for \allowbreak$\phi[\top,\dots,\top]$ over $2^N$ such that $\tau(x)\subsetneq 2^N$ for all leaf nodes $x$ labeled with $ r_i$ ($1\leq i\leq m$) in $\mathfrak{T}_\phi$.

By condition (iii),  there exist $1\leq j_1\leq\dots\leq j_k\leq  m$ such that (\ref{non_uniform_def_eq4}) holds for the set $N$. On the other hand, for each $j_a$ ($1\leq a\leq m$), as $2^N\nsubseteq\tau(x)$ holds for every leaf node $x$ labeled with $ r_{j_a}$, we have that  $\tau(x)\models\Theta^\star_{2^N}$, i.e.,
\(\tau(x)\models \mathsf{f}(x)[\alpha_1,\dots,\alpha_m],\)
where each $\alpha_i$ is defined as in Equation (\ref{non_uniform_def_eq3}). Thus, by Lemma \ref{truth_function_determine_Xi}, $\tau$ is also a truth function for $\phi[\alpha_1,\dots,\alpha_m]$
over $2^N$,  thereby  
\(2^N\models \phi [\alpha_1,\dots,\alpha_m].\)
Thus, by (\ref{non_uniform_def_eq2}), we obtain that $2^N\models \cnt(\alpha_1,\dots,\alpha_m)$,
which contradicts (\ref{non_uniform_def_eq4}).
\end{proof}

Finally, we are in a position to derive the main result of the paper as a corollary of the above theorem.

\begin{thm}\label{bor_non_uniform_def}
Intuitionistic disjunction is not uniformly definable in \PD.
\end{thm}
\begin{proof}
It suffices to check that intuitionistic disjunction satisfies conditions (i)-(iii)  of Theorem \ref{non_uniform_def}. Condition (i) is satisfied, since, e.g., $\bot\bor \top\not\models \bot$ and $\top\bor \bot\not\models\bot$. Condition (ii) is satisfied since, e.g., $\models\top\bor\top$.  Lastly, for any finite set $N$ of propositional variables, $2^N\not\models\Theta^\star_{2^N}\bor\Theta^\star_{2^N}$, giving condition (iii).
\end{proof}

We have already proved that intuitionistic implication is not uniformly definable in \PD in Theorem \ref{imp_non_uniform_def} by observing that intuitionistic implication is not monotone. In fact, the non-uniform definability of intuitionistic implication in \PD also follows from Theorem \ref{non_uniform_def}, as intuitionistic implication also satisfies conditions (i)-(iii).
 Indeed, we have that (i) $ \bot\to \bot\not\models \bot$, (ii) $\models\top\to\top$ and (iii) $2^N\not\models\top\to\Theta^\star_{2^N}$. 

 Finally, we summarize the results obtained in this section as a corollary concerning compositional translatability between \Inql and \PD. One may compare this corollary with Corollary \ref{pid_non_uniform_def}.
 
\begin{cor}
$\Inql\leq\PD$, whereas $\Inql\not\leqc\PD$.
\end{cor}
\begin{proof}
By Theorems \ref{max_logic}, \ref{imp_non_uniform_def} and \ref{bor_non_uniform_def}.
\end{proof}


\section{Concluding remarks}


 Team semantics was originally devised (in the context of independence-friendly logic) by \cite{Hodges1997a,Hodges1997b} to meet one of the fundamental needs of logic and language, namely ``\emph{compositionality}'' (see, e.g., \cite{compositionality_J97,compositionality_Hodges01,PaginWesterstahl2010} for an overview).  However, the result of this paper, as well as those in \cite{ivano_msc}, \cite{Pietro_uniform} show that there is a distinction between definability and uniform definability, and between translatability and compositional translatability in team semantics. This phenomenon seems to indicate that the compositionality or uniformity on another level is lost in team semantics. Although it is commented in \cite{Hodges2016} that the results presented in this paper do not actually have conflicts with the  notion of compositionality given in \cite{Hodges2012}, in the author's opinion, there are yet a lot more to be clarified regarding this subtle issue.
 
We finish by mentioning that it is an open problem whether, on the other hand, \PD is compositionally translatable into \Inql. The dependence atoms are uniformly definable in \Inql, since $\dep(p_1,\dots,p_k,q)\equiv (p_1\vee\neg p_1)\wedge\dots\wedge (p_k\vee\neg p_k)\to (q\vee\neg q)$. But  whether the tensor $\sor$ is uniformly definable in \Inql is open. We conjecture that it is not, and note that the argument in this paper does not seem to work for the logic \Inql, as contexts for \Inql, especially those that contain intuitionistic implication, are not in general monotone in the sense of Lemma \ref{context_pd_monotone} (cf. the proof of Theorem \ref{imp_non_uniform_def}).

\bibliographystyle{rsl}

%

\clearpage

\end{document}